\documentclass[11pt]{article}
\usepackage{graphicx}
\usepackage{amsmath}
\usepackage{amsfonts}
\usepackage{amssymb}

\def\sqr#1#2{{\vcenter{\vbox{\hrule height.#2pt
              \hbox{\vrule width.#2pt height#1pt \kern#1pt \vrule width.#2pt}
              \hrule height.#2pt}}}}
\def\signed #1{{\unskip\nobreak\hfil\penalty50
              \hskip2em\hbox{}\nobreak\hfil#1
              \parfillskip=0pt \finalhyphendemerits=0 \par}}
\def\endpf{\signed {$\sqr69$}}

\def\ns{\noalign{\medskip}}

\newtheorem{lemma}{Lemma}[section]
\newtheorem{theorem}{Theorem}[section]

\newtheorem{proposition}{Proposition}[section]

\newtheorem{definition}{Definition}[section]

\newtheorem{remark}{Remark}[section]
\newtheorem{example}{Example}[section]

\makeatletter
   
   \@addtoreset{equation}{section}
\makeatother
\date{}

\begin{document}

\title{{\bf Classifications of Linear Controlled
Systems}\thanks{This work was partially supported by the NSF of
China under grants 10525105 and 10771149. The author gratefully
acknowledges Professor Xu Zhang for his help.\medskip}}

\author{Jing Li
\thanks{School of Mathematics,  Sichuan University,
Chengdu 610064, China. {\small\it E-mail:} {\small\tt
lijing4924@163.com}.}}
\date{\quad}
\maketitle

\begin{abstract}
\noindent This paper is devoted to a study of linear, differential
and topological classifications for linear controlled systems
governed by ordinary differential equations. The necessary and
sufficient conditions for the linear and topological equivalence are
given. It is also shown that the differential equivalence is the
same as the linear equivalence for the linear controlled systems.
\end{abstract}

\bigskip


\newpage

\newpage

{\begin{footnotesize}

\tableofcontents

\end{footnotesize}}

\section{Introduction}

In this paper, we study the classification problems for the
following time-invariant linear controlled system
\begin{equation}\label{11}\overset{\cdot }{x}\left(
t\right)=Ax(t)+Bu(t),\quad t\geq 0,
\end{equation}
where $\overset{\cdot }{x}=\frac{dx}{dt}$, $x(t)\in\mathbb{R}^n$ is
the {\it state variable}, $u(t)\in\mathbb{R}^m$ is the {\it control
variable}, and $A$ and $B$ are real matrices of dimensions $n\times
n$ and $n\times m$, respectively. The admissible control set is
chosen to be $L^1_{\text{loc}}(\mathbb{R}^+;\mathbb{R}^m)$. Since
(\ref{11}) is uniquely determined by the pair of matrices $A$ and
$B$, we denote it simply by $(A,B)$.

Classification is a basic problem in science. Elements in the same
equivalent class have some similar properties (such as
controllability and the number of the efficient controls of
controlled systems in our case). Thus, in order to know the
properties of all elements in the same class, we only need to study
some special one (for example, the system with the canonical form
that will be introduced later).

Consider now the following two systems
\begin{equation}\label{12}\overset{\cdot }{x}\left(
t\right)=A_1x(t)+B_1u(t)
\end{equation}
and
\begin{equation}\label{13}\overset{\cdot }{y}(t)=A_2y(t)+B_2v(t).
\end{equation}
Here, $x(t)$, $y(t)\in\mathbb{R}^n$ are state variables, $u(t)$,
$v(t)\in\mathbb{R}^m$ are control variables, and $A_i$, $B_i$ ($i=1,
2$) are real matrices of dimensions $n\times n$ and $n\times m$,
respectively. We introduce the following:

\begin{definition}\label{d11}
1) Systems (\ref{12}) and (\ref{13}) are called topologically
equivalent if there exists a (vector-valued) function $F(x,u)\equiv
(H(x,u),G(x,u))$, where $H(x,u)\in
C(\mathbb{R}^n\times\mathbb{R}^m;\mathbb{R}^n)$ and $G(x,u)\in
C(\mathbb{R}^n\times\mathbb{R}^m;\mathbb{R}^m)$, such that

\begin{enumerate}
\item[i)]$F(\cdot,\cdot)$ is a homeomorphism from $\mathbb{R}^n\times\mathbb{R}^m$ to $\mathbb{R}^n\times\mathbb{R}^m$ (henceforth we denote the inverse function of $F(x,u)$ by
$F^{-1}(y,v)\equiv (Z(y,v),W(y,v))$);

\item[ii)]The transformation $(y(t),v(t))=(H(x(t),u(t)),G(x(t),u(t)))$ brings
(\ref{12}) to (\ref{13}), and the transformation
$(x(t),u(t))=(Z(y(t),v(t)),W(y(t),v(t)))$ brings (\ref{13}) to
(\ref{12}).
\end{enumerate}

\noindent 2) Systems (\ref{12}) and (\ref{13}) are called
differentially equivalent if $F(x,u)$ is an $C^1$ diffeomorphism
from $\mathbb{R}^n\times\mathbb{R}^m$ to
$\mathbb{R}^n\times\mathbb{R}^m$.

\noindent 3) Systems (\ref{12}) and (\ref{13}) are called linearly
equivalent if $F(x,u)$ is a linear isomorphism from
$\mathbb{R}^n\times\mathbb{R}^m$ to
$\mathbb{R}^n\times\mathbb{R}^m$.
\end{definition}

Several remarks are in order.

\begin{remark}\label{r11}
The transformation $(y(t),v(t))=(H(x(t),u(t)),G(x(t),u(t)))$ brings
(\ref{12}) to (\ref{13}) means that: if $x(t)$ is the solution of
(\ref{12}) with the initial datum $x(0)=x^0$ and the control
$u(\cdot)\in C([0,+\infty);\mathbb{R}^m)$, then by the
transformation $F(\cdot,\cdot)$, $y(t)=H(x(t),u(t))$ is the solution
of (\ref{13}) with the initial datum $y(0)=H(x(0),u(0))=H(x^0,u(0))$
and the control $v(t)=G(x(t),u(t))$.
\end{remark}

\begin{remark}\label{r12}
It is easy to check that the linear, differential and topological
equivalence are actually equivalence relations, by which we mean, as
usual, they are symmetric, reflexive and transitive. In the sequel,
we call $(H(x,u),G(x,u))$ the equivalence transformation from system
(\ref{12}) to system (\ref{13}).
\end{remark}

\begin{remark}\label{r13}
It is clear that the following relation holds:
\begin{eqnarray*}
\text{Linear equivalence} \Rightarrow \text{Differential
equivalence} \Rightarrow \text{Topological equivalence}.
\end{eqnarray*}
\end{remark}

\begin{remark}\label{r14}
Since any ordinary differential equation can be regarded as a
controlled system without effective control (i.e. $B=0$ in
(\ref{11})), the classification results in this paper also apply to
ordinary differential equations (ODEs for short). For autonomous
linear ODEs, our result coincides with those in \cite{SZ}.
\end{remark}

The classification problem for the linear completely controllable
system was studied by Brunovsky \cite{Brunovsky}. In
\cite{Brunovsky}, Brunovsky introduced the concept of feedback
equivalence and showed that there are only finitely many feedback
equivalence classes and each of which can be represented by a simple
canonical form. Later, Liang (\cite{Liang}) extended this result to
the general controlled linear system.

There exist extensive works on the classification of controlled
systems by means of {\it linear} equivalence transformations (cf.
\cite{BFZ,Brunovsky,Liang,LNP,SZ1}). As for the {\it nonlinear}
transformations, we refer to \cite{Elkin} and \cite{TR} for some
results on nonlinear controlled system $\overset{\cdot
}{x}=f(x)+g(x)u$ with $C^{\infty}$ smooth transformation
$y=\phi(x)$, $v=\alpha(x)+\beta(x)u$. In this paper, we will adopt
the methods introduced in \cite{SZ} (which is addressed to the
classifications of ordinary differential equations) and our
transformations are neither necessarily linear nor need high
regularity.

We claim that the function $H(x,u)$ in Definition \ref{d11} has the
following fundamental property.

\begin{proposition}\label{p11}
If $F(x,u)=(H(x,u),G(x,u))$ is a topological equivalence
transformation from system (\ref{12}) to system (\ref{13}), then
$H(x,u)$ is independent of $u$ and we can simply write it as $H(x)$.
Furthermore, $H(x)$ is a homeomorphism from $\mathbb{R}^n$ to
$\mathbb{R}^n$.
\end{proposition}

We refer the reader to Appendix A for a proof of Proposition
\ref{p11}. It follows from Proposition \ref{p11} that Definition
\ref{d11} can be reduced to the following simpler one.

\begin{definition}\label{d12}
1) Systems (\ref{12}) and (\ref{13}) are called topologically
equivalent if there exists a (vector-valued) function $F(x,u)\equiv
(H(x),G(x,u))$, where $H(x)\in C(\mathbb{R}^n;\mathbb{R}^n)$ and
$G(x,u)\in C(\mathbb{R}^n\times\mathbb{R}^m;\mathbb{R}^m)$, such
that

\begin{enumerate}
\item[i)]$F(\cdot,\cdot)$ is a homeomorphism from $\mathbb{R}^n\times\mathbb{R}^m$ to $\mathbb{R}^n\times\mathbb{R}^m$ (henceforth we denote the inverse function of $F(x,u)$ by
$F^{-1}(y,v)\equiv (H^{-1}(y),W(y,v))$);

\item[ii)]The transformation $(y(t),v(t))=(H(x(t)),G(x(t),u(t)))$ brings
(\ref{12}) to (\ref{13}), and the transformation
$(x(t),u(t))=(H^{-1}(y(t)),W(y(t),v(t)))$ brings (\ref{13}) to
(\ref{12}).
\end{enumerate}

\noindent 2) Systems (\ref{12}) and (\ref{13}) are called
differentially equivalent if $F(x,u)$ is an $C^1$ diffeomorphism
from $\mathbb{R}^n\times\mathbb{R}^m$ to
$\mathbb{R}^n\times\mathbb{R}^m$.

\noindent 3) Systems (\ref{12}) and (\ref{13}) are called linearly
equivalent if $F(x,u)$ is a linear isomorphism from
$\mathbb{R}^n\times\mathbb{R}^m$ to
$\mathbb{R}^n\times\mathbb{R}^m$.
\end{definition}

In what follows, we discuss the equivalence classes of system
(\ref{11}) in the sense of Definition \ref{d12}. The main results of
this paper are follows. First, it is shown that the linear
equivalence given in Definition \ref{d12} coincides with the
feedback equivalence which will be recalled later. Next, we study
the differential classification of system (\ref{11}). Since the
differential equivalence transformation $F(x,u)$ is an $C^1$
diffeomorphism, we can find a linear isomorphism with the aid of the
derivative of $F$. It turns out that the differential classification
is the same as the linear classification. Finally, we address to the
topological classification of system (\ref{11}). Technically, the
discussion is reduced to the classification for a completely
controllable system and an ODE, respectively. The crucial point is
to classify linear completely controllable systems in our sense of
topological equivalence. The difficulty consists in the fact that
very little is known about the properties of the equivalence
transformation. Thus, it is hard to find a homeomorphism. To
overcome this, we study the canonical form and use the contradiction
argument to avoid constructing the equivalence transformation
directly. Our result read: once two linear completely controllable
systems are topologically equivalent, then they are linearly
equivalent.

The rest of this paper is organized as follows. Some preliminary
knowledge are recalled in Section \ref{sp}. In Section \ref{sl}, we
discuss the linear classification for system (\ref{11}).  Section
\ref{sd} is devoted to analyzing the differential classification for
system (\ref{11}). In Section \ref{st}, we study the topological
classification and give a necessary and sufficient condition for
this classification. Section \ref{st1} and Appendix A are devoted
respectively to prove Proposition \ref{p45} and Proposition
\ref{p11}.

\section{Some preliminaries}\label{sp}

In this section, we present some preliminary results, which will
play a role in the sequel.

\subsection{Feedback classification of linear controlled systems}

In \cite{Liang}, Liang extended the concept of feedback equivalence
which was introduced by Brunovsky for linear completely controllable
systems to general controlled linear systems.

\begin{definition}\label{d2.3.1}
Systems (\ref{12}) and (\ref{13}) are called feedback equivalent if
there exist matrices $O$, $Q$ and $L$ of dimensions $n\times n$,
$m\times m$ and $m\times n$, respectively, with $O$ and $Q$ being
nonsingular, such that
\begin{eqnarray*}
 A_2=O^{-1}A_1O+O^{-1}B_1L, \quad \quad B_2=O^{-1}B_1Q.
\end{eqnarray*}
\end{definition}

Put $k=\text{rank}\left (B,AB,\cdots,A^{n-1}B \right)$. Clearly,
$k\leq n$. It is well-known that system (\ref{11}) is completely
controllable if and only if $k=n$.

Following Brunovsky (\cite{Brunovsky}), we introduce two sequences
for linear controlled systems. Put $r_0=\text{rank}B$ and
$r_j=\text{rank}\left (B,AB,\cdots,A^jB \right)-\text{rank}
\left(B,AB,\cdots,A^{j-1}B \right)$ for $ 1\leq j\leq n-1$. Define
$R(A,B)\equiv \{ r_j \}_{j=0}^{n-1}$. Denote by $L_j$ the linear
subspace of $\mathbb{R}^n$ spanned by the column vectors of
$\left(B, AB, \cdots, A^jB \right)$. Denote by $A_j$ the orthogonal
complement of $L_{j-1}$ in $L_j$ and by $\pi_j(b)$ the orthogonal
projection of a vector $b\in L_j$ into $A_j$. One can choose k
linear independent vectors from the column vectors of $\left
(B,AB,\cdots,A^{n-1}B \right)$ to construct a set $\cal S$, such
that the vectors $\left \{ \pi_j(A^jb_i);\;A^jb_i\in {\cal S},\ j \
\text{fixed} \right \}$ span $A_j$ and if $A^jb_i\not\in {\cal S}$,
then $A^{j+1}b_i\not\in {\cal S}$, where $b_i\in \mathbb{R}^{n}$ is
the $i$-th column of $B$. Associate every column $b_i$ with a number
$p_i$, such that $A^jb_i\in {\cal S}$ for $0\leq j\leq p_i-1$, but
$A^{p_i}b_i\not\in {\cal S}$. By re-ordering suitably the columns of
$B$, one can achieve that $p_1\geq p_2\geq\cdots\geq p_m$. Define
$P(A,B)\equiv\{ p_i \}_{i=1}^{m}$.

\medskip

It is easy to deduce from Lemma 1 in Section 2 of \cite{Brunovsky}
that

\begin{lemma}\label{l21}
The finite sequences $R(A,B)= \{ r_j \}_{j=0}^{n-1}$ and $P(A,B)=\{
p_i \}_{i=1}^{m}$ have the following properties:

\begin{enumerate}
\item[1)] $0\leq r_j\leq m$,\quad  $r_0\geq r_1\geq \cdots  \geq
r_{p_1-1}>0$,\quad $r_j=0$ for $j\geq p_1$,\quad $\displaystyle
\sum_{j=0}^{n-1}r_j=k$;

\item[2)] $0\leq p_i\leq n$,\quad $p_1\geq p_2\geq \cdots  \geq
p_{r_0}>0$,\quad $p_i=0$ for $i> r_0$,\quad $\displaystyle
\sum_{i=1}^mp_i=k$;

\item[3)] $P(A_1,B_1)=P(A_2,B_2)$ if and only if
$R(A_1,B_1)=R(A_2,B_2)$.
\end{enumerate}
\end{lemma}

\begin{lemma}\label{l44}
{\rm (Theorem 1 in Section 2 of \cite{Brunovsky})} Assume that
systems (\ref{12}) and (\ref{13}) are completely controllable. Then
they are feedback equivalent if and only if $R(A_1,B_1)=R(A_2,B_2)$
(or equivalently $P(A_1,B_1)=P(A_2,B_2)$).
\end{lemma}

We derive from \cite{Liang} the following result.

\begin{theorem}\label{t22}
Any system $(A,B)$ is feedback equivalent to a system
$(\widetilde{A},\widetilde{B})$ of the form:
\begin{eqnarray*}
\widetilde{A}= \left[\begin{array}{cc}
 C & 0\\
 0 & M
\end{array}\right],
\quad\quad \widetilde{B}= \left[\begin{array}{c}
 D\\
 0
\end{array}\right],
 \end{eqnarray*}
where $M$ is an $(n-k)\times (n-k)$ Jordan matrix of the form
(recall $k=\text{rank}\left (B,AB,\cdots,A^{n-1}B \right)$)
\begin{eqnarray}\label{22}
&& M= \left[
  \begin{array}{ccc}
M^- & 0 & 0\\
0 & M^+ & 0\\
0 & 0 & M^0
\end{array}
\right],
\end{eqnarray}
for which the real parts of eigenvalues of $M^-$, $M^+$ and $M^0$
are negative, positive and zero, respectively; $C$ and $D$ are
respectively $k\times k$ and $k\times m$ matrices of the forms
\begin{eqnarray*}
 C= \left[
  \begin{array}{cccc}
J_{p_1} & 0 & \cdots & 0\\
0 & J_{p_2} & \cdots & 0\\
\vdots & \vdots & \ddots & \vdots\\
0 & 0 &\cdots & J_{p_{r_0}}
\end{array}
\right], \quad \quad D= \left[
  \begin{array}{ccccccc}
e_{p_1} & 0  & \cdots & 0 & 0 & \cdots & 0\\
0 & e_{p_2} & \cdots & 0 & 0 & \cdots & 0\\
\vdots & \vdots &\ddots & \vdots & \vdots && \vdots\\
0 & 0 &\cdots & e_{p_{r_0}} & 0 & \cdots & 0
\end{array}
\right],
\end{eqnarray*}
where $\{p_i\}_{i=1}^m=P(A,B)$, $J_q$ and $e_q$ are $q \times q$ and
$q \times 1$ matrices, respectively
\begin{eqnarray*}
J_q= \left[
  \begin{array}{ccccc}
0 & 1 & 0 & \cdots & 0\\
0 & 0 & 1 & \cdots & 0\\
\vdots & \vdots & \vdots & \ddots & \vdots\\
0 & 0 & 0 & \cdots & 1\\
0 & 0 & 0 & \cdots & 0
\end{array}
\right],\quad \quad
 e_q=
 \left[
  \begin{array}{c}
0 \\
\vdots\\
0\\
1
\end{array}
\right].
\end{eqnarray*}
\end{theorem}

\begin{remark}\label{r22}
In the sequel, we call the pair $(\widetilde{A},\widetilde{B})$ in
Theorem \ref{t22} the canonical form of $(A,B)$.
\end{remark}

\subsection{Classification of autonomous linear ordinary differential
equations}

Consider the following two autonomous linear ordinary differential
equations
\begin{equation}\label{2.1.3}\overset{\cdot }{x}\left(
t\right)=A_1x(t)
\end{equation}
and
\begin{equation}\label{2.1.4}\overset{\cdot }{y}(t)=A_2y(t),
\end{equation}
where $x(t)$, $y(t) \in \mathbb {R}^n$, and $A_1$, $A_2$ are real
matrices of dimensions $n\times n$.

\begin{definition}\label{d2.1.1}
{\rm(Definition 6.1 in Chapter 2 of \cite[p. 26]{SZ})} Systems
(\ref{2.1.3}) and (\ref{2.1.4}) are called topologically equivalent
if there exists a (vector-valued) function $H(x)\in
C(\mathbb{R}^n;\mathbb{R}^n)$ such that the transformation
$y(t)=H(x(t))$ brings (\ref{2.1.3}) to (\ref{2.1.4}).
\end{definition}

Denote by $n^-$, $n^+$ and $n^0$ respectively the numbers of matrix
$M$'s eigenvalues with negative real parts, positive real parts and
zero real parts (similar notations $n_i^-$, $n_i^+$ and $n_i^0$ are
for matrix $M_i$ ($i=1,2$)).

\begin{lemma}\label{l43}
{\rm(Theorem 11.1 in Chapter 3 of \cite[p. 49]{SZ})} ODEs $(M_1,0)$
and $(M_2,0)$ are topologically equivalent if and only if
$(n^-_1,n^+_1,n^0_1)=(n^-_2,n^+_2,n^0_2)$ and the matrices $M_1^0$
and $M_2^0$ are similar. Here, $M_i$ ($i=1,2$) are matrices in the
form (\ref{22}).
\end{lemma}

\subsection{Lebesgue covering theorem}

For any topological space, the Lebesgue covering dimension is
defined to be $n$ if $n$ is the smallest integer for which the
following holds: any open cover has a refinement (a second cover
where each element is a subset of an element in the first cover)
such that no point is included in more than $n + 1$ elements. If no
such $n$ exists, then the dimension is infinite.

\begin{theorem}\label{t23} {\rm(Lebesgue covering theorem, c.f. \cite{HW})} The Lebesgue covering dimension
coincides with the affine dimension of a finite simplicial complex.
\end{theorem}

Theorem \ref{t23} indicates that the dimension do not change via a
homeomorphism transformation.

\section{Linear classification}\label{sl}

This section is addressed to the linear classification of system
(\ref{11}). The main result in this section is stated as follows.

\begin{theorem}\label{t21}
Systems (\ref{12}) and (\ref{13}) are linearly equivalent if and
only if there are matrices $O$, $Q$ and $L$ of dimensions $n\times
n$, $m\times m$ and $m\times n$, respectively, with $O$ and $Q$
being nonsingular, such that
\begin{eqnarray}\label{21}
 A_2=O^{-1}A_1O+O^{-1}B_1L, \quad \quad B_2=O^{-1}B_1Q.
\end{eqnarray}
\end{theorem}

{\it Proof.} {\it Sufficiency.} It is easy to check that
\begin{eqnarray*}
\left[\begin{array}{c}
 y\\
 v
\end{array}\right]
= \left[\begin{array}{c}
 H(x)\\
 G(x,u)
\end{array}\right]
= \left[\begin{array}{cc}
 O^{-1} & 0\\
 -Q^{-1}LO^{-1} & Q^{-1}
\end{array}\right]
\left[\begin{array}{c}
 x\\
 u
\end{array}\right]
\end{eqnarray*}
is the equivalence transformation from system (\ref{12}) to system
(\ref{13}). Therefore these two systems are linearly equivalent.

\medskip

{\it Necessity.} Assume systems (\ref{12}) and (\ref{13}) are
linearly equivalent and denote the equivalence transformation by
\begin{eqnarray*}
 \left[\begin{array}{c}
 y\\
 v
\end{array}\right]
= \left[\begin{array}{c}
 H(x)\\
 G(x,u)
\end{array}\right]
=\left[\begin{array}{cc}
 \ell_1 & 0\\
 \ell_2 & \ell_3
\end{array}\right]
\left[\begin{array}{c}
 x\\
 u
\end{array}\right],
\end{eqnarray*}
where $\left[\begin{array}{cc}
 \ell_1 & 0\\
 \ell_2 & \ell_3
\end{array}\right]$ is a nonsingular $(n+m)\times(n+m)$ constant matrix. Hence $\ell_1$ and
$\ell_3$ are nonsingular. Set
$$ O=\ell_1^{-1},\quad
Q=\ell_3^{-1},\quad L=-\ell_3^{-1}\ell_2\ell_1^{-1},$$ then $A_i$,
$B_i$ (i=1, 2) satisfy relation $(\ref{21})$.
\endpf

\begin{remark}\label{r21}
Theorem \ref{t21} shows that linear equivalence defined in this
paper is actually the same as feedback equivalence introduced in
\cite{Brunovsky} (recall Definition \ref{d2.3.1} for feedback
equivalence).

The novelty in our Definition \ref{d12} for the linear equivalence
is that we discuss the equivalence by the property of transformation
and define three types of classifications in a unified way.
\end{remark}

\section{Differential classification}\label{sd}

This section is devoted to the differential classification of system
(\ref{11}). The main result in this section is as follows.

\begin{theorem}
Systems (\ref{12}) and (\ref{13}) are differentially equivalent if
and only if they are linearly equivalent.
\end{theorem}

{\it Proof.} By Remark \ref{r13}, it suffices to show the `` only if
'' part.

Denote by $(H(x),G(x,u))$ the differential equivalence
transformation from system (\ref{12}) to system (\ref{13}). The
solutions of systems (\ref{12}) and (\ref{13}) can be expressed
respectively as
\begin{eqnarray*}
 x(t)=e^{A_1t}x^0+\int_0^t e^{A_1(t-s)}B_1u(s)ds
\end{eqnarray*}
and
\begin{eqnarray}\label{32}
y(t)= H(x(t))=e^{A_2t}y^0+\int_0^t e^{A_2(t-s)}B_2G(x(s),u(s))ds.
\end{eqnarray}
Taking $t=0$ in (\ref{32}) gives $ H\left (x^0 \right)=y^0. $ Hence,
(\ref{32}) can be rewritten as
\begin{eqnarray}\label{33}
 H(x(t))=e^{A_2t}H \left(x^0 \right)+\int_0^t e^{A_2(t-s)}B_2G(x(s),u(s))ds.
\end{eqnarray}
Setting $u\equiv 0$, differentiating (\ref{33}) with respect to
$x^0$, and then setting $x^0=0$, we arrive at
\begin{eqnarray}\label{34}
 D_xH(0)e^{A_1t}=e^{A_2t}D_xH(0)+\int_0^t e^{A_2(t-s)}B_2D_xG(0,0)e^{A_1s}ds.
\end{eqnarray}
Furthermore, differentiating (\ref{34}) with respect to $t$ and
setting $t=0$, we get
\begin{eqnarray}\label{35}
 D_xH(0)A_1=A_2D_xH(0)+B_2D_xG(0,0).
\end{eqnarray}
Similarly, setting $x^0=0$ and $u\equiv u^0$, differentiate
(\ref{33}) with respect to $t$, and then setting $t=0$, it follows
\begin{eqnarray*}
D_xH(0)B_1u^0=A_2H(0)+B_2G\left(0,u^0\right).
\end{eqnarray*}
Differentiating the above formula with respect to $u^0$ and setting
$u^0=0$, we find
\begin{eqnarray}\label{39}
D_xH(0)B_1=B_2D_uG(0,0).
\end{eqnarray}

Since $F(x,u)$ is a diffeomorphism from $\mathbb{R}^n \times
\mathbb{R}^m$ to itself, we see that $ F^{-1}(F(x,u))=(x,u)$.
Consequently,
\begin{eqnarray*}
D_{(y,v)}F^{-1}(F(x,u)) D_{(x,u)}F(x,u)=I,
\end{eqnarray*}
which implies
\begin{eqnarray*}
\text{det}D_{(y,v)}F^{-1}(F(x,u))\text{det}D_{(x,u)}F(x,u)=1.
\end{eqnarray*}
From this formula, we deduce that
\begin{eqnarray*}
D_{(x,u)}F(0,0) = \left[\begin{array}{cc}
 D_xH(0) & 0 \\
 D_xG(0,0) & D_uG(0,0)
\end{array}\right]
\end{eqnarray*}
is nonsingular. Thus, $D_xH(0)$ and $D_uG(0,0)$ are nonsingular.
Therefore, (\ref{35}) and (\ref{39}) yield
\begin{eqnarray}\label{37}
 A_1=D_xH(0)^{-1}A_2D_xH(0)+D_xH(0)^{-1}B_2D_xG(0,0),
\end{eqnarray}
and
\begin{eqnarray}\label{318}
B_1=D_xH(0)^{-1}B_2D_uG(0,0).
\end{eqnarray}
Combining $(\ref{37})$ and $(\ref{318})$,  and noting Theorem 2.1,
we arrive at the desired result. This completes the proof of Theorem
3.1.
\endpf

\begin{remark}
Theorem 3.1 indicates that for system (\ref{11}), its differential
classification  is the same as its linear classification. However,
as we will show later, its linear classification is different from
its topological classification. Indeed, Example \ref{e43} in Section
\ref{st} gives two systems, which are not linearly equivalent but
topologically equivalent.
\end{remark}

\section{Topological classification}\label{st}

This section is devoted to the topological classification of system
(\ref{11}).

To begin with, let us introduce some notations. Put
$$\displaystyle
k_i=\text{rank}\left(B_i,A_iB_i,\cdots,A^{n-1}_iB_i\right),\quad
i=1, 2. $$
In virtue of Theorem \ref{t22}, systems $(A_i,B_i)$ are linearly
equivalent to their canonical forms
\begin{eqnarray}\label{451}
 (\widetilde{A}_i,\widetilde{B}_i)=\left(
\begin{array}{cc}
 \left[
 \begin{array}{cc}
   C_i & 0\\
   0 & M_i\\
 \end{array}
 \right],&
 \left[
 \begin{array}{c}
   D_i\\
   0
 \end{array}
 \right]
\end{array}
\right), \quad i=1, 2,
\end{eqnarray}
with $C_i$, $D_i$ and $M_i$ being $k_i\times k_i$, $k_i\times m$ and
$(n-k_i)\times (n-k_i)$ matrices, respectively.
Denote $r=\text{rank}B_1$, $s=\text{rank}B_2$,
$\{p_i\}_{i=1}^m=P(A_1,B_1)$ and $\{q_i\}_{i=1}^m=P(A_2,B_2)$ Then
$C_i$ and $D_i$ can be expressed as follows (recall that $J_q$ and
$e_q$ are defined in Theorem \ref{t22})
\begin{eqnarray}\label{452}
 C_1= \left[
  \begin{array}{cccc}
J_{p_1} & 0 & \cdots & 0\\
0 & J_{p_2} & \cdots & 0\\
\vdots & \vdots & \ddots & \vdots\\
0 & 0 &\cdots & J_{p_{r}}
\end{array}
\right], \quad\quad D_1= \left[
  \begin{array}{ccccccc}
e_{p_1} & 0  & \cdots & 0 & 0 & \cdots & 0\\
0 & e_{p_2} & \cdots & 0 & 0 & \cdots & 0\\
\vdots & \vdots &\ddots & \vdots & \vdots && \vdots\\
0 & 0 &\cdots & e_{p_r} & 0 & \cdots & 0
\end{array}
\right],
\end{eqnarray}
and
\begin{eqnarray}\label{453}
 C_2= \left[
  \begin{array}{cccc}
J_{q_1} & 0 & \cdots & 0\\
0 & J_{q_2} & \cdots & 0\\
\vdots & \vdots & \ddots & \vdots\\
0 & 0 &\cdots & J_{q_s}
\end{array}
\right], \quad\quad D_2= \left[
  \begin{array}{ccccccc}
e_{q_1} & 0  & \cdots & 0 & 0 & \cdots & 0\\
0 & e_{q_2} & \cdots & 0 & 0 & \cdots & 0\\
\vdots & \vdots &\ddots & \vdots & \vdots && \vdots\\
0 & 0 &\cdots & e_{q_s} & 0 & \cdots & 0
\end{array}
\right].
\end{eqnarray}

\begin{remark}\label{r41}
For any topological equivalence transformation $(H(x),G(x,u))$ from
system (\ref{12}) to system (\ref{13}), without loss of generality,
we can always assume $H(0)=0$ in the following discussion.
Otherwise, one may choose $(H(x)-H(0),G(x,u)-G(0,0))$ to be the new
topological equivalence transformation.
\end{remark}

The main result in this section is stated as follows.

\begin{theorem}\label{t41}
Assume that the canonical forms of systems (\ref{12}) and (\ref{13})
are respectively
\begin{eqnarray*}
(\widetilde{A}_i,\widetilde{B}_i)=
 \left(
\begin{array}{cc}
 \left[
 \begin{array}{cc}
   C_i & 0 \\
   0 & M_i\\
 \end{array}
 \right],&
 \left[
 \begin{array}{c}
   D_i\\
   0
 \end{array}
 \right]
\end{array}
\right), \quad i=1, 2, \end{eqnarray*} given by (\ref{451}). Then,
systems (\ref{12}) and (\ref{13}) are topologically equivalent if
and only if
\begin{eqnarray}\label{411}
&k_1= k_2,\\ \ns\label{412} &R(A_1,B_1)=R(A_2,B_2)\ \quad \text{(or
equivalently \  $P(A_1,B_1)=P(A_2,B_2)$)},\\ \ns\label{413}
&(n^-_1,n^+_1,n^0_1)=(n^-_2,n^+_2,n^0_2),
\end{eqnarray} and the matrices $M_1^0$ and $M_2^0$ are similar.
\end{theorem}

For better orientation, we divide the proof of this theorem into
several propositions which concern with the necessary conditions for
the topological equivalence.

\begin{proposition}\label{p41}
If systems (\ref{12}) and (\ref{13}) are topologically equivalent,
then $k_1=k_2$.
\end{proposition}

{\it Proof.} We use the contradiction argument and assume that
$k_1>k_2$.

By Theorem \ref{t22}, systems $(A_i,B_i)$ are linearly equivalent to
systems $(\widetilde{A}_i,\widetilde{B}_i)$ $(i=1, 2)$. Hence, by
Remark \ref{r13}, without loss of generality, we assume that systems
$(A_i,B_i)$ $(i=1, 2)$ are in the canonical form (\ref{451}).

Denote by $(H(x),G(x,u))$ the topological equivalence transformation
from (\ref{12}) to (\ref{13}). Let $x(0)=0$ and fix a $t_0>0$. Since
system (\ref{12}) is composed by a completely controllable system
$ (C_1,D_1) $ and an ODE $ (M_1,0) $, it is easy to see that
\begin{eqnarray*}
\Theta_1:=\left\{x\left(t_0;0,u(\cdot)\right);\; u(\cdot)\in
C([0,+\infty);\mathbb{R}^m)\right\} =\mathbb {R}^{k_1}\times
\{0\}\subset \mathbb {R}^n.
\end{eqnarray*}
Similarly, for any initial datum $y^0\in\mathbb{R}^n$, one finds
\begin{eqnarray*}\begin{array}{ll}
\displaystyle\Theta_2:=\left\{y(t_0;y^0,v(\cdot));\;  v(\cdot)\in
C\left([0,+\infty);\mathbb{R}^m\right)\right\}\\ \ns\quad\ \
=\mathbb{R}^{k_2}\times \left \{\text{a fixed point in}\
\mathbb{R}^{n-k_2}
 \right\}\subset\mathbb{R}^n.\end{array}
\end{eqnarray*}

Thanks to Theorem \ref{t23}, we see that $H(\Theta_1)$ is a
$k_1$-dimensional topological manifold since $H(x)$ is a
homeomorphism from $\mathbb{R}^n$ to $\mathbb{R}^n$. On the other
hand, since transformation $(y(t),v(t))=(H(x(t),G(x(t),u(t)))$
brings (\ref{12}) to (\ref{13}), one has $H(\Theta_1)\subset
\Theta_2$. Noting the dimension of $\Theta_2$ is $k_2$, we conclude
that the dimension of $H(\Theta_1)$ is at most $k_2$. Noting however
that $k_1>k_2$, we arrive at a contradiction. This completes the
proof of Proposition \ref{p41}.
\endpf

\begin{proposition}\label{p42*}
Let $E_i$, $K_i$ and $N_i$ be real matrices of dimensions $k\times
k$, $k\times m$ and $(n-k)\times (n-k)$,  respectively. Assume that
system
\begin{equation}\label{441} \left[\begin{array}{c}
   \overset{\cdot }{x}_1(t)\\
   \overset{\cdot }{x}_2(t)
 \end{array}
 \right]
=  \left[
 \begin{array}{cc}
   E_1 & 0\\
   0 & N_1\\
 \end{array}
 \right]
\left[\begin{array}{c}
   x_1(t)\\
   x_2(t)
 \end{array}
 \right]
+\left[\begin{array}{c}
   K_1\\
   0
 \end{array}
 \right]u(t)
\end{equation}
is topologically equivalent to system
\begin{equation}\label{442} \left[\begin{array}{c}
   \overset{\cdot }{y}_1(t)\\
   \overset{\cdot }{y}_2(t)
 \end{array}
 \right]
=  \left[
 \begin{array}{cc}
   E_2 & 0\\
   0 & N_2\\
 \end{array}
 \right]
\left[\begin{array}{c}
   y_1(t)\\
   y_2(t)
 \end{array}
 \right]
+\left[\begin{array}{c}
   K_2\\
   0
 \end{array}
 \right]v(t),
\end{equation}
and the topological equivalence transformation from $(\ref{441})$ to
$(\ref{442})$ is
\begin{eqnarray*}\left\{\begin{array}{ll}
y_1=h_1(x_1,x_2),\\ \ns y_2=h_2(x_1,x_2),\\ \ns v=G(x_1,x_2,u).
\end{array}\right.
\end{eqnarray*}
And assume systems $\displaystyle\overset{\cdot
}{x}_1(t)=E_1x_1(t)+K_1u(t) $ and $ \displaystyle\overset{\cdot
}{y}_1(t)=E_2y_1(t)+K_2v(t) $ are completely controllable. Then,
$y_2=h_2(x_1,x_2)$ is independent of $x_1$ and $h_2(x_2)$ is a
homeomorphism from $\mathbb{R}^{n-k}$ to $\mathbb{R}^{n-k}$.
\end{proposition}

{\it Proof.} The proof is divided into several steps.

\medskip

{\it Step 1.} We show that $h_2(x_1,x_2)$ is independent of $x_1$,
that is
$$h_2(x_1,x_2)\equiv h_2(x_2),\quad \forall x_1\in\mathbb{R}^k,\ \forall
x_2\in\mathbb{R}^{n-k}.$$
In fact, assuming this is not the case, we deduce that there exist
$a^0,b^0\in\mathbb{R}^k$ and $\bar{x}_2\in\mathbb{R}^{n-k}$ such
that
\begin{equation}\label{447} h_2 \left (a^0, \bar{x}_2\right) \neq h_2 \left(b^0,
\bar{x}_2\right).\end{equation}

We use $\left(x_1 (t), x_2 (t)\right)$ to denote the solution of
(\ref{441}) with initial datum $\displaystyle\left(x_1^0,
x_2^0\right)$ and control $u(t)$. Denote by $\left(\hat{x}_1 (t),
\hat{x}_2 (t)\right)$ the solution of (\ref{441}) with initial datum
$\displaystyle\left(\hat{x}_1^0, \hat{x}_2^0\right)$ and control
$\hat{u}(t)$. Similar notations are for system (\ref{442}).

Fix a $t_0 >0$ and a $x_1^0\in\mathbb{R}^k$. Taking
$\hat{x}_1^0=x_1^0$ and $\hat{x}_2^0=x_2^0 = e^{-N_1t_0}\bar{x}_2$,
one finds
\begin{equation}\label{448}\hat{x}_2(t_0)=e^{N_1t_0}\hat{x}_2^0=\bar{x}_2=e^{N_1t_0}x_2^0=x_2(t_0).\end{equation}
Since system $\overset{\cdot }{x}_1(t)=E_1x_1(t)+K_1u(t)$ is
completely controllable, there exist two controls $u(t)$ and
$\hat{u}(t)$ such that $x_1(t_0) = a^0$ and $\hat{x}_1 (t_0) = b^0$.
Hence, combining (\ref{447}) and (\ref{448}), it follows
\begin{equation}\label{449}
\hat{y}_2(t_0)=h_2\left(\hat{x}_1(t_0),\hat{x}_2(t_0)\right)=h_2\left(b^0,\bar{x}_2\right)\neq
h_2\left(a^0,\bar{x}_2\right)=
h_2\left(x_1(t_0),x_2(t_0)\right)=y_2(t_0).
\end{equation}

However, by
$$\hat{y}_2^0=h_2\left(\hat{x}_1^0,\hat{x}_2^0\right)=h_2\left(x_1^0,x_2^0\right)=y_2^0,$$
we obtain
$$\hat{y}_2(t_0)=e^{N_2t_0}\hat{y}_2^0=e^{N_2t_0}y_2^0=y_2(t_0),$$
which contradicts (\ref{449}). Therefore, $h_2 (x_1, x_2)$ is
independent of $x_1$ and henceforth we simply denote  it by $h_2
(x_2)$.

\medskip

{\it Step 2.} We now show that $h_2(x_2)$ is a homeomorphism from
$\mathbb{R}^{n-k}$ to $\mathbb{R}^{n-k}$.

First of all, we prove that $h_2$ is invertible. Indeed, we denote
the inverse function of \linebreak $(h_1(x_1,x_2),h_2(x_2))$ by
$(z_1(y_1,y_2),z_2(y_1,y_2))$ (recall that, by Proposition
\ref{p11}, $(h_1(x_1,x_2),h_2(x_2))$ is invertible). Similar to Step
1, one finds
$$z_2(y_1,y_2)\equiv z_2(y_2),\quad \forall y_1\in\mathbb{R}^k,\ \forall
y_2\in\mathbb{R}^{n-k}.$$ Therefore,
$$(x_1,x_2)=(z_1(h_1(x_1,x_2),h_2(x_2)),z_2(h_2(x_2))),\quad \forall x_1\in\mathbb{R}^k,\ \forall x_2\in\mathbb{R}^{n-k},$$
which yields $$x_2=z_2(h_2(x_2)),\quad \forall
x_2\in\mathbb{R}^{n-k}.$$ Similarly, we have
$$y_2=h_2(z_2(y_2)),\quad \forall y_2\in\mathbb{R}^{n-k}.$$ Hence $h_2^{-1}=z_2$.

Next, noting $h_2$ and $z_2$ are continuous functions, we conclude
that $h_2(x_2)$ is a homeomorphism from $\mathbb{R}^{n-k}$.\endpf

\begin{proposition}\label{p42}
Under the assumption of Proposition \ref{p42*}, systems
\begin{equation}\label{443}
\displaystyle\overset{\cdot }{x}_1(t)=E_1x_1(t)+K_1u(t)
\end{equation}
and
\begin{equation}\label{444}
\displaystyle\overset{\cdot }{y}_1(t)=E_2y_1(t)+K_2v(t)
\end{equation}
are topologically equivalent, and ODEs
\begin{equation}\label{445}
\displaystyle\overset{\cdot }{x}_2(t)=N_1x_2(t)
\end{equation}
and
\begin{equation}\label{446}
\displaystyle\overset{\cdot }{y}_2(t)=N_2y_2(t)
\end{equation}
are topologically equivalent as well.
\end{proposition}

{\it Proof.} It is easy to see that the transformation
$y_2(t)=h_2(x_2(t))$ brings (\ref{445}) to (\ref{446}).
Consequently,  $y_2=h_2(x_2)$ is an equivalence transformation from
system (\ref{445}) to system (\ref{446}), and therefore these two
ODEs are topologically equivalent.

We show that system (\ref{443}) is topologically equivalent to
system (\ref{444}).

From Proposition \ref{p42*}, the inverse function of
$(h_1(x_1,x_2),h_2(x_2);G(x_1,x_2,u))$ is of the form\linebreak
$(z_1(y_1,y_2),z_2(y_2);W(y_1,y_2,v))$.

First of all, we claim that $(h_1(x_1, 0);G(x_1, 0, u))$ is a
homeomorphism from $\mathbb{R}^k \times \mathbb{R}^m$ to
$\mathbb{R}^k \times \mathbb{R}^m$ and its inverse function is
$(z_1(y_1,0);W(y_1,0,v))$. Indeed, Remark \ref{r41} implies $ h_2
(0)=0$. Since $h_2(x_2)$ is one-to-one, $x_2 =0$ if and only if $y_2
=0$. Noting $(h_1(x_1,x_2),h_2(x_2);G(x_1,x_2,u))$ is a
homeomorphism of $\mathbb {R}^n\times\mathbb {R}^m$, we see that
$(h_1(x_1,0),h_2(0) ;G(x_1,0,u))$ is a homeomorphism of $\mathbb
{R}^k \times \{ 0 \} \times \mathbb {R}^m$ and its inverse function
is $(z_1(y_1, 0), z_2 (0);W(y_1, 0, v))$. Therefore, $(h_1(x_1,
0);G(x_1, 0, u))$ is a homeomorphism of $\mathbb{R}^k \times
\mathbb{R}^m$ and its inverse function is $(z_1(y_1,0);W(y_1,0,v))$.

Next, for any $(x_1(t),u(t))$ satisfying system (\ref{443}), it is
easy to see that $(x_1(t),0,u(t))$ satisfies system (\ref{441}).
Since the transformation $(h_1(x_1,x_2),h_2(x_2);G(x_1,x_2,u))$
brings (\ref{441}) to (\ref{442}), we conclude that
$(h_1(x_1(t),0),h_2(0),G(x_1(t),0,u(t))$ satisfies (\ref{442}).
Hence $(h_1(x_1(t),0),\linebreak G(x_1(t),0,u(t)))$ satisfies
(\ref{444}), which yields that the transformation
$(h_1(x_1,0);G(x_1,0,u))$ brings (\ref{443}) to (\ref{444}).

Similarly, the inverse transformation $(z_1(y_1,0);W(y_1,0,v))$ of
$(h_1(x_1,0);G(x_1,0,u))$ brings (\ref{444}) to (\ref{443}).

Thus,  $(y_1,v)=(h_1(x_1,0);G(x_1,0,u))$ is an equivalence
transformation from system (\ref{443}) to system (\ref{444}), and
therefore these two controlled systems are topologically
equivalent.
\endpf

\begin{remark}
Theorem \ref{t22} indicates that in the sense of linear equivalence,
system (\ref{11}) can be reduced to a completely controllable system
$\overset{\cdot }{\xi}(t)=C\xi(t)+Du(t)$ and an ODE $ \overset{\cdot
}{\eta}(t)=M\eta(t)$ (i.e. a controlled system without effective
control).

Moreover, by Proposition \ref{p42}, it suffices to study the
topological classification of the completely controllable system and
the ODE, respectively.
\end{remark}

Next, we will discuss some other necessary conditions for the
topological equivalent with the aid of the following two lemmas.

\begin{lemma}\label{l41}
If $(H(x),G(x,u))$ is a topological equivalence transformation from
system $(\widetilde{A}_1,\widetilde{B}_1)$ to system
$(\widetilde{A}_2,\widetilde{B}_2)$, and
$$
\bar{x}= \Big(\overset{p_1 \text{ times}}{\overbrace{\bar{x}_1, 0,
\cdots, 0}},\  \overset{p_2 \text{ times}}{\overbrace{\bar{x}_{p_1
+1}, 0, \cdots, 0}},\  \cdots, \ \overset{p_r \text{
times}}{\overbrace{\bar{x}_{\sum_{j=1}^{r-1} p_j +1}, 0, \cdots,
0}},\ \overset{n-k_1 \text{
times}}{\overbrace{0,\cdots,0}}\Big)^{\top},
$$
then the corresponding point $\bar{y} = H\left(\bar{x}\right)$ has
the following form
\begin{equation}\label{471}
\bar{y}= \Big(\overset{q_1 \text{ times}}{\overbrace{\bar{y}_1, 0,
\cdots, 0}},\ \overset{q_2 \text{ times}}{\overbrace{ \bar{y}_{q_1
+1}, 0, \cdots, 0}},\ \cdots,\ \overset{q_s \text{
times}}{\overbrace{ \bar{y}_{\sum_{j=1}^{s-1} q_j +1}, 0, \cdots,
0}},\ \overset{n-k_2 \text{
times}}{\overbrace{0,\cdots,0}}\Big)^{\top}.
\end{equation}
Here, $x^{\top}$ denotes the transpose of a vector (or matrix) $x$
and $(\widetilde{A}_i,\widetilde{B}_i)$ ($i=1, 2$) are canonical
forms given by (\ref{451}).
\end{lemma}

{\it Proof.} Noting the special form of matrices $\widetilde{A}_1$
and $\widetilde{B}_1$, the solution of system
$(\widetilde{A}_1,\widetilde{B}_1)$ associated to $x(0)=\bar{x}$ and
$ u(t)\equiv 0$ is
$$x(t)\equiv \bar{x},\quad\forall t\geq 0.$$
Therefore, the corresponding $y(t)=H(x(t))\equiv
H\left(\bar{x}\right)$ and $v(t)=G(x(t),u(t))\equiv G(\bar{x},0)$
are two constants. We denote
$\bar{y}=(\bar{y}_1,\cdots,\bar{y}_n)^{\top} =
H\left(\bar{x}\right)=(h_1(\bar{x}),\cdots,h_n(\bar{x}))^{\top}$ and
$ v(t)\equiv (v_1^0,\cdots,v_m^0)^{\top}$.

\smallskip

Substituting $ \xi(t)\equiv(\bar{y}_1,\cdots,\bar{y}_{\sum_{j=1}^{s}
q_j })^{\top}$ and $ v(t)\equiv (v_1^0,\cdots,v_m^0)^{\top}$ into
system $\overset{\cdot }{\xi}(t)=C_2\xi(t)+D_2v(t)$, we deduce that
\begin{eqnarray*}
(0,\cdots,0)^{\top}=\Big(\overset{q_1 \text{
times}}{\overbrace{\bar{y}_2,\cdots,\bar{y}_{q_1},v_1^0}},\
\overset{q_2 \text{
times}}{\overbrace{\bar{y}_{q_1+2},\cdots,\bar{y}_{q_1+q_2},v_2^0}},\
\cdots,\ \overset{q_s \text{ times}}{\overbrace{
\bar{y}_{\sum_{j=1}^{s-1} q_j +2},\cdots,\bar{y}_{\sum_{j=1}^{s} q_j
},v_s^0}}\Big)^{\top},
\end{eqnarray*}
which gives
\begin{equation}\label{472}
\left(\bar{y}_1,\cdots,\bar{y}_{\sum_{j=1}^{s} q_j }\right)^{\top}=
\Big(\overset{q_1 \text{ times}}{\overbrace{\bar{y}_1, 0, \cdots,
0}},\ \overset{q_2 \text{ times}}{\overbrace{ \bar{y}_{q_1 +1}, 0,
\cdots, 0}},\  \cdots, \ \overset{q_s \text{
times}}{\overbrace{\bar{y}_{\sum_{j=1}^{s-1} q_j +1}, 0, \cdots,
0}}\Big)^{\top}.
\end{equation}

By Proposition \ref{p42*}, we see that $(\bar{y}_{\sum_{j=1}^{s}
q_j+1 },\cdots,\bar{y}_n)^{\top}$ is independent of
$(\bar{x}_1,\cdots,\bar{x}_{\sum_{j=1}^{r} p_j })^{\top}$. Hence
$$\left(\bar{y}_{\sum_{j=1}^{s} q_j+1 },\cdots,\bar{y}_n\right)^{\top}
=\left(h_{\sum_{j=1}^{s} q_j+1
}(\bar{x}),\cdots,h_n(\bar{x})\right)^{\top}=\left(h_{\sum_{j=1}^{s}
q_j+1 }(0),\cdots,h_n(0)\right)^{\top}.$$
Since Remark \ref{r41} implies $H(0)=0$,  the above formula yields
\begin{equation}\label{473}\left(\bar{y}_{\sum_{j=1}^{s} q_j+1
},\cdots,\bar{y}_n\right)^{\top}
=\left(0,\cdots,0\right)^{\top}.\end{equation}

Combining (\ref{472}) and (\ref{473}) we arrive at the desired
(\ref{471}). This completes the proof of Lemma \ref{l41}.
\endpf

\begin{lemma}\label{l42}
Assume that $(H(x),G(x,u))$ is a topological equivalence
transformation from $(\widetilde{A}_1,\widetilde{B}_1)$ to
$(\widetilde{A}_2,\widetilde{B}_2)$, where
$(\widetilde{A}_i,\widetilde{B}_i)$ ($i=1, 2$) are canonical forms
given by (\ref{451}). Suppose $x(0) = 0$ and $u(t)$ is continuous.
Then
\begin{enumerate}
\item[1)] the solution $x(t)=(x_1(t),\cdots,x_n(t))^{\top}$ and the control
$u(t)= (u_1(t),\cdots,u_m(t))^{\top}$ of system
$(\widetilde{A}_1,\widetilde{B}_1)$ satisfy
$$
u_{i}  = \frac{d x_{\sum_{j=1}^{i} p_j} }{dt} = \frac{d^2
x_{\sum_{j=1}^{i} p_j -1} }{dt^2} = \cdots = \frac{d^{p_{i}}
x_{\sum_{j=1}^{{i}-1} p_j -p_{i} +1}}{dt^{p_{i}}}, \quad {i} = 1,
\cdots, r;
$$

\item[2)] the corresponding solution $ y(t) = H(x(t))$ and the corresponding
control $v(t) = G(x(t), u(t))$ of system
$(\widetilde{A}_2,\widetilde{B}_2)$ satisfy
$$
v_{i}  =\frac{d y_{\sum_{j=1}^{i} q_j}}{dt} = \frac{d^2
y_{\sum_{j=1}^{i} q_j -1}}{dt^2} = \cdots = \frac{d ^{q_{i}}
y_{\sum_{j=1}^{i} q_j -q_{i} +1}}{dt^{q_{i}}}, \quad {i} = 1,
\cdots, s.
$$
\end{enumerate}
\end{lemma}

{\it Proof.} From $x(0 )=0$, noting the special structure of
matrices $\widetilde{A}_1$ and $\widetilde{B}_1$, we have
\begin{eqnarray*}
\begin{array}{l}
\displaystyle x_{\sum_{j=1}^{i} p_j} (t) = \int_0^t u_{i} (\tau) d
{\tau}, \quad x_{\sum_{j=1}^{i} p_j -1} (t) = \int_0^t
\int_0^{\tau} u_{i} (\varsigma) d {\varsigma} d{\tau}, \cdots,\\
\ns\displaystyle x_{\sum_{j=1}^{i} p_j -p_{i} +1} (t) = \int_0^t
\int_0^{\tau} \cdots \int_0^{\varrho} u_{i} (\nu) d {\nu}  \cdots
d{\varsigma}d {\tau}, \quad {i} =1, \cdots, r,
\end{array}
\end{eqnarray*}
which yields the first conclusion.

\medskip

By Remark \ref{r41} we have $y(0)=H(x(0))=H(0)=0$. Similar to the
above, noting $v(t)= G(x(t),u(t))$ is continuous, we see that
\begin{eqnarray*}
\begin{array}{l}
\displaystyle y_{\sum_{j=1}^{i} q_j} (t) = \int_0^t v_{i} (\tau) d
{\tau}, \quad y_{\sum_{j=1}^{i} q_j -1} (t) = \int_0^t
\int_0^{\tau} v_{i} (\varsigma) d {\varsigma} d{\tau}, \cdots,\\
\ns \displaystyle y_{\sum_{j=1}^{i} q_j -q_{i} +1} (t) = \int_0^t
\int_0^{\tau} \cdots \int_0^{\varrho} v_{i} (\sigma) d {\sigma}
\cdots d{\varsigma}d {\tau}, \quad {i} =1, \cdots, s,
\end{array}
\end{eqnarray*}
which gives the second conclusion.
\endpf

\begin{proposition}\label{p43}
If systems (\ref{12}) and (\ref{13}) are topologically equivalent,
then $\text{rank}B_1=\text{rank}B_2$.
\end{proposition}

{\it Proof.} We use the contradiction argument. If
$\text{rank}B_1\neq\text{rank}B_2$, we denote $\text{rank}B_1=r$,
$\text{rank}B_2=s$ and assume $r>s$.

By Theorem \ref{t22} and Remark \ref{r13}, without loss of
generality, we assume systems $(A_i,B_i)$ $(i=1, 2)$ are in the
canonical form (\ref{451}).

Denote by $(H(x),G(x,u))$ the equivalence transformation from
(\ref{12}) to (\ref{13}). Set
\begin{eqnarray*}\begin{array}{ll}
\Theta:=\Big\{\Big(x_1^0, 0, \cdots, 0; x_{p_1 +1}^0, 0, \cdots, 0;
\cdots; x_{\sum_{j=1}^{r-1} p_j +1}^0, 0, \cdots,
0;0,\cdots,0\Big)^{\top} ; \\ \ns\displaystyle\qquad\qquad x_1^0,\
x_{p_1+1}^0\ ,\cdots, x_{\sum_{j=1}^{r-1} p_j
+1}^0\in\mathbb{R}\Big\}.\end{array}
\end{eqnarray*}
Thanks to Lemma \ref{l41}, we have
\begin{equation}\label{474}
H(\Theta)=\Big\{\Big(y_1^0, 0, \cdots, 0; y_{q_1 +1}^0, 0, \cdots,
0; \cdots; y_{\sum_{j=1}^{s-1} q_j +1}^0, 0, \cdots,
0;0,\cdots,0\Big)^{\top}\Big\}.
\end{equation}

By Theorem \ref{t23}, $H(\Theta)$ is a $r$-dimensional topological
manifold since $\Theta$ is a $r$-dimensional topological manifold.
On the other hand, (\ref{474}) gives $ y^0$ is at most a
$s$-dimensional topological manifold, which contradicts the fact
$r>s$. This completes the proof of Proposition \ref{p43}.
\endpf

\begin{remark}\label{r42}
Proposition \ref{p43} guarantees that the numbers of the effective
controls of two equivalence systems coincide. Therefore, an ODE
(\text{rank}$B_1=0$) can not be topologically equivalent to any
controlled system with effective control (\text{rank}$B_2\neq 0$).
\end{remark}

\begin{proposition}\label{p45}
If systems $\overset{\cdot }{x}(t)=C_1x(t)+D_1u(t)$ and
$\overset{\cdot }{y}(t)=C_2y(t)+D_2v(t)$ are topologically
equivalent, then they are linearly equivalent. Here, $C_i$ and $D_i$
($i=1, 2$) are $k\times k$ and $k\times m$ matrices given by
(\ref{452}) and (\ref{453}), respectively.
\end{proposition}

The proof of Proposition \ref{p45} will be given in Section
\ref{st1}.

\begin{remark}
It is a direct consequence Proposition \ref{p45} that if systems
(\ref{12}) and (\ref{13}) are completely controllable, then they are
topologically equivalent if and only if they are linearly
equivalent. However, as shown in Example \ref{e43}, this is not
always the case for general controlled systems.
\end{remark}

{\it Proof of Theorem \ref{t41}.} Without loss of generality, we
assume
\begin{eqnarray*}
 A_i= \left[
 \begin{array}{cc}
   C_i & 0\\
   0 & M_i\\
 \end{array}
 \right], \quad\quad B_i= \left[
 \begin{array}{c}
   D_i\\
   0
 \end{array}
 \right],\quad i=1, 2.
\end{eqnarray*}

{\it Sufficiency.} By (\ref{411}) and (\ref{412}), the completely
controllable systems $(C_1,D_1)$ and $ (C_2,D_2)$ are linearly
equivalent follows directly from Lemma \ref{l44}.

By (\ref{413}) and noticing that $M_1^0$ is similar to $M_2^0$, from
Lemma \ref{l43} we see that ODEs $(M_1,0)$ and $(M_2,0)$ are
topologically equivalent.

Therefore, system (\ref{12}) is topologically equivalent to system
(\ref{13}).

\medskip

{\it Necessity.} If systems (\ref{12}) and (\ref{13}) are
topologically equivalent, then Proposition \ref{p41} guarantees that
(\ref{411}) holds. Furthermore, Proposition \ref{p42} implies that
the completely controllable systems $(C_1,D_1)$ and $ (C_2,D_2)$ are
topologically equivalent, and ODEs $(M_1,0)$ and $(M_2,0)$ are
topologically equivalent. By Proposition \ref{p45} and Lemma
\ref{l44}, it is easy to see that condition (\ref{412}) holds.
Thanks to Lemma \ref{l43}, we get that (\ref{413}) holds and $M_1^0$
is similar to $M_2^0$.
\endpf

\bigskip

Theorem \ref{t41} shows that the ODE part brings difference between
the topological classification and the linear classification. Let us
consider a example which support this point of view.

\begin{example}\label{e43}
Systems \begin{eqnarray} \label{41}\displaystyle\overset{\cdot
}{x}(t)= \left[
 \begin{array}{cc}
   0 & 0\\
   0 & a
 \end{array}
 \right]x(t)
+\left[
 \begin{array}{c}
   1\\
   0
 \end{array}
 \right]u(t)
\end{eqnarray}
and
\begin{eqnarray} \label{42}\displaystyle\overset{\cdot }{y}(t)=
\left[
 \begin{array}{cc}
   0 & 0\\
   0 & -1
 \end{array}
 \right]y(t)
+\left[
 \begin{array}{c}
   1\\
   0
 \end{array}
 \right]v(t)
\end{eqnarray} are topologically equivalent, but not linearly equivalent.
Here $a<0$ and $a\neq -1$.
\end{example}

{\it Proof.} By Theorem \ref{t41}, it is easy to check that these
two systems are topologically equivalent.

Now, we prove that systems (\ref{41}) and (\ref{42}) are not
linearly equivalent. Otherwise, due to Theorem \ref{t21}, there
exist two nonsingular matrices $O= \left[\begin{array}{cc}o_1&
o_2\\o_3& o_4
\end{array}\right]^{-1}$, $Q=[q]$, and a matrix $L=[l_1,l_2]$ such
that
\begin{equation}\label{414}
\left[\begin{array}{cc}
   0 & 0\\
   0 & -1
\end{array}\right]
= \left[\begin{array}{cc}
 o_1& o_2\\o_3& o_4
\end{array}\right]
\left[\begin{array}{cc}
   0 & 0\\
   0 & a
\end{array}\right]
\left[\begin{array}{cc}
 o_1& o_2\\o_3& o_4
\end{array}\right]^{-1}
+ \left[\begin{array}{cc}
 o_1& o_2\\o_3& o_4
\end{array}\right]
\left[\begin{array}{c}
   1\\
   0
\end{array}\right]
\left[\begin{array}{cc}
 l_1& l_2
\end{array}\right]
\end{equation}
and
\begin{equation}\label{415} \left[\begin{array}{c}
   1\\
   0
\end{array}\right]
= \left[\begin{array}{cc}
 o_1& o_2\\o_3& o_4
\end{array}\right]
\left[\begin{array}{c}
   1\\
   0
\end{array}\right]
[q].
\end{equation}
From (\ref{415}), we get $o_3=0$. Thus, the identity (\ref{414})
gives
\begin{eqnarray*}
&& \left[\begin{array}{cc}
   0 & 0\\
   0 & -1
\end{array}\right]
= \left[\begin{array}{cc}
 o_1l_1& \frac{ao_2}{o_4}+o_1l_2\\\displaystyle 0 & a
\end{array}\right],
\end{eqnarray*}
which contradicts the fact $a\neq -1$. Hence (\ref{41}) and
(\ref{42}) are not linearly equivalent.
\endpf

\section{Proof of Proposition \ref{p45}}\label{st1}

This section is addressed to a proof of Proposition \ref{p45}. To
begin with, we consider the following lemma, which is a special case
of Proposition \ref{p45} with $\text{rank}D_1=\text{rank}D_2=2$.

\begin{lemma}\label{p44}
If $p_1 + p_2 = q_1 + q_2$ and $ p_2 < q_2$, then the two completely
controllable systems
\begin{eqnarray} \label{416***}\displaystyle\overset{\cdot }{x}(t)=
\left[\begin{array}{cc} J_{p_1} & 0\\0&J_{p_2}
\end{array}\right]x(t)
+\left[\begin{array}{cc}e_{p_1}&0\\0&e_{p_2}\end{array}\right]u(t)
\end{eqnarray}
and
\begin{eqnarray}\label{417***}\displaystyle\overset{\cdot }{y}(t)=
\left[\begin{array}{cc} J_{q_1} & 0\\0&J_{q_2}
\end{array}\right]y(t)
+\left[\begin{array}{cc}e_{q_1}&0\\0&e_{q_2}\end{array}\right]v(t)
\end{eqnarray}
are not topologically equivalent (recall that $J_q$ and $e_q$ are
defined in Theorem \ref{t22}).
\end{lemma}

Clearly, this lemma is a special case of Proposition \ref{p45}.

{\it Proof.} We use the contradiction argument and divide the proof
into several steps.

\medskip

{\it Step 1.}  Assume that (\ref{416***}) is topologically
equivalent to (\ref{417***}) and the equivalence transformation is
$(H(x),G(x,u))$.

For any $u_2(t)\in C([0,+\infty);\mathbb{R}^1)$, the solutions of
system (\ref{416***}) associated to $x(0) = 0$ and $u(t) = (0, u_2
(t))^{\top}$ is
\begin{eqnarray}\label{418***}
\begin{array}{ll} x(t)& =\Big(\overset{p_1 \text{
times}}{\overbrace{0, \cdots, 0}},\ \overset{p_2 \text{
times}}{\overbrace{ x_{p_1 +1} (t), \cdots, x_{p_1 + p_2}
(t)}}\Big)^{\top}
\\ \ns & = \displaystyle \Big(\overset{p_1 \text{
times}}{\overbrace{0,\cdots,0}},\ \overset{p_2 \text{
times}}{\overbrace{\int_0^t \cdots \int_0^{\varsigma} u_2(\varrho) d
{\varrho} \cdots d{\tau},\cdots, \int_0^t u_2(\tau) d\tau}}
\Big)^{\top},\quad\forall t\geq 0
\end{array}
\end{eqnarray}
Denote by $\displaystyle H(x(t)) = y(t) = (y_1 (t), \cdots, y_{p_1 +
p_2} (t))^{\top}$ the corresponding solution of system
(\ref{417***}) via the transformation $\displaystyle H(x) = (h_1(x),
\cdots, h_{p_1 + p_2}(x))^{\top}$. It is obvious that $y(t)$ only
depends on $ x_{p_1 +1} (t),  x_{p_2 +2} (t), \cdots, x_{p_1 + p_2}
(t) $. Define
\begin{equation*}
\Phi_i (\lambda_1, \cdots, \lambda_{p_2}) = h_{q_1+i} (\overset{p_1
\text{ times}}{\overbrace{0, \cdots, 0}},\ \overset{p_2 \text{
times}}{\overbrace{ \lambda_1, \cdots, \lambda_{p_2} }}), \quad i=1,
\cdots, q_2.
\end{equation*}
Clearly, $\Phi_i$ $(i=1, \cdots, q_2)$ are continuous functions in
$(\lambda_1, \cdots, \lambda_{p_2}) \in \mathbb{R}^{p_2}$ and
\begin{equation}\label{418.5**}
y_{q_1 + i} (t) = \Phi_i \left(x_{p_1 +1} (t), x_{p_1 +2} (t),
\cdots, x_{p_1 + p_2}(t)\right), \quad i=1, \cdots, q_2.
\end{equation}

\medskip

{\it Step 2.} We claim that: for any $(\lambda_1, \cdots,
\lambda_{p_2}) \in \mathbb{R}^{p_2}$, it holds
\begin{equation}\label{419***}
\Phi_{p_2} (\lambda_1,  \cdots, \lambda_{p_2 -1}, \lambda_{p_2})
\equiv \Phi_{p_2} (\lambda_1,  \cdots, \lambda_{p_2 -1}, 0).
\end{equation}
In fact, assuming this is not the case, we conclude that there exist
$a_1, \cdots, a_{p_2 -1}\in\mathbb{R}$ and $ a_{p_2} \in
\mathbb{R}\backslash\{ 0\}$ such that
\begin{equation*}
\Phi_{p_2} (a_1, \cdots, a_{p_2 -1}, a_{p_2}) \neq \Phi_{p_2} (a_1,
\cdots, a_{p_2 -1}, 0).
\end{equation*}
Denote $b = \Phi_{p_2} (a_1, \cdots, a_{p_2 -1}, a_{p_2}) -
\Phi_{p_2} (a_1, \cdots, a_{p_2 -1}, 0)$. Since $\Phi_{p_2}$ is
continuous, there exists a constant $d > 0$ such that for any
 $(\lambda_1,
\cdots, \lambda_{p_2 -1}) \in [a_1 -d, a_1 +d] \times \cdots \times
[a_{p_2 -1} -d, a_{p_2 -1} +d] $, we have
\begin{equation}\label{419.5***}
\left| \Phi_{p_2} (\lambda_1, \cdots, \lambda_{p_2 -1}, a_{p_2}) -
\Phi_{p_2} (a_1, \cdots, a_{p_2 -1}, a_{p_2}) \right| \le
\frac{|b|}{2}.
\end{equation}
Since $\Phi_{p_2 +1}$ is continuous and $[a_1 -d, a_1 +d] \times
\cdots \times [a_{p_2 -1} -d, a_{p_2 -1} +d] \times [-|a_{p_2}|,
|a_{p_2}|] \subset \mathbb{R}^{p_2}$ is a compact set, there exists
a constant $c > 0$ such that for any $(\lambda_1, \cdots,
\lambda_{p_2 -1}, \lambda_{p_2}) \in [a_1 -d, a_1 +d] \times \cdots
\times [a_{p_2 -1} -d, a_{p_2 -1} +d] \times [-|a_{p_2}|,
|a_{p_2}|]$, it holds
\begin{equation}\label{420***}
\left |\Phi_{p_2 +1}(\lambda_1, \cdots, \lambda_{p_2 -1},
\lambda_{p_2})\right | \le c.
\end{equation}

Put $$t_0 = \min\left \{2, 1+ \displaystyle \frac{d}{|a_2| + |a_3| +
\cdots + |a_{p_2}|}, 1 + \frac{|b|}{4c}\right \}.$$
Choose a $C^{\infty}$ function $\varphi$ satisfying that: for
$\displaystyle0\leq t\leq \frac{1}{4}$, $\varphi(t)\equiv 0$, and
for $\displaystyle t\geq \frac{1}{2}$,
\begin{equation}\label{420.5***}
\varphi (t) = a_1 + a_2 (t-1) + \frac{a_3}{2!} (t-1)^2 + \cdots +
\frac{a_{p_2 -1}}{(p_2 -2)!} (t-1)^{p_2 -2} + \frac{a_{p_2}
(t-1)^{p_2}}{p_2 !(t_0-1)}.
\end{equation}
Setting $u_2 (t) = \varphi^{(p_2)} (t)$ and recalling $x(0)=0$, by
Lemma \ref{l42} and (\ref{418.5**}), we obtain that for
$\displaystyle t\geq\frac{3}{4}$,
\begin{eqnarray}\label{421***}
\begin{array}{ll}
x_{p_1 +1} (t) = \varphi (t),\\ \ns
\displaystyle x_{p_1 +2} (t) =\overset{\cdot}{\varphi}(t)\\
\ns\displaystyle = a_2 + a_3 (t-1) + \cdots + \frac{a_{p_2 -1}}{(p_2
-3)!} (t-1)^{p_2 -3} +
\frac{a_{p_2} (t-1)^{p_2 -1}}{(p_2 -1) !(t_0-1)}, \\ \ns \cdots \\
\ns
\displaystyle x_{p_1 + p_2 -1} (t) = \varphi^{(p_2-2)} (t) = a_{p_2
-1} + \frac{a_{p_2} (t-1)^2}{2!(t_0-1)}, \\ \ns
\displaystyle x_{p_1 + p_2} (t)= \varphi^{(p_2-1)} (t) =
\frac{a_{p_2} (t-1)}{t_0-1};
\end{array}
\end{eqnarray}
and
\begin{eqnarray}\label{425***}
\begin{array}{ll}
\displaystyle \Phi_{p_2 +1} \left(\varphi (t),
\overset{\cdot}{\varphi} (t), \cdots, \varphi^{(p_2 -1)} (t)\right)
= y_{q_1 + p_2 + 1} (t) =\overset{\cdot} {y}_{q_1 + p_2} (t)
\\ \ns
\displaystyle = \frac{d \Phi_{p_2} \left(\varphi (t),
\overset{\cdot}{\varphi} (t), \cdots, \varphi^{(p_2 -1)}
(t)\right)}{d t}.
\end{array}
\end{eqnarray}
From (\ref{420.5***}), (\ref{421***}) and the choice of $t_0$ we get
that for $1 \le t \le t_0$
\begin{eqnarray*}
\begin{array}{ll}
\displaystyle \left| \varphi (t) - a_1 \right| = \left| a_2 (t-1) +
\frac{a_3}{2!} (t-1)^2 + \cdots + \frac{a_{p_2 -1}}{(p_2 -2)!}
(t-1)^{p_2 -2} + \frac{a_{p_2} (t-1)^{p_2}}{p_2 !(t_0-1)} \right| \\
\ns
\displaystyle\qquad\qquad\ \ \le  |a_2| (t - 1) + |a_3| (t-1) +
\cdots + |a_{p_2 -1}| (t-1) + \frac{ |a_{p_2}| (t-1)^2}{t_0 -1} \\
\ns
\displaystyle\qquad\qquad\ \  \le  (|a_2| + |a_3| + \cdots +
|a_{p_2}|) (t_0 - 1)\leq d,
\\ \ns  \cdots \\ \ns
\displaystyle \left | \varphi^{(p_2-2)} (t) - a_{p_2 -1} \right |  =
\left| \frac{a_{p_2} (t-1)^2}{2!(t_0-1)} \right| \le |a_{p_2}| (t_0
- 1) \leq d, \\ \ns
\displaystyle \left |\varphi^{(p_2-1)} (t) \right |  = \left|
\frac{a_{p_2} (t-1)}{t_0-1} \right| \le |a_{p_2}|,
\end{array}
\end{eqnarray*}
which gives that
\begin{eqnarray}\label{420.6***}
\begin{array}{c}
|  \varphi (t) - a_1  | \le d,\ \cdots,\ | \varphi^{(p_2-2)} (t) -
a_{p_2 -1} | \le d, \quad  | \varphi^{(p_2-1)} (t) | \le |a_{p_2}|,
\quad \forall 1 \le t \le t_0.
\end{array}
\end{eqnarray}
Combining (\ref{420***}) and (\ref{420.6***}), one has
\begin{eqnarray}\label{420.7***}
\displaystyle \left |\Phi_{p_2 +1} \left(\varphi (t),
 \cdots, \varphi^{(p_2 -2)} (t), \varphi^{(p_2 -1)} (t)\right) \right |
\le c, \quad \forall 1 \le t \le t_0.
\end{eqnarray}

However, by (\ref{419.5***}) and (\ref{420.6***}), and setting
$t=t_0$, we arrive at
\begin{equation}\label{420.8***}
\left | \Phi_{p_2}\left (\varphi (t_0), \cdots, \varphi^{(p_2 -2)}
(t_0),a_{p_2}\right) - \Phi_{p_2}(a_1,  \cdots, a_{p_2-1}, a_{p_2})
\right| \le \frac{|b|}{2}.
\end{equation}
Recalling $b = \Phi_{p_2} (a_1, \cdots, a_{p_2 -1}, a_{p_2}) -
\Phi_{p_2} (a_1, \cdots, a_{p_2 -1}, 0)$, from (\ref{420.5***}) and
(\ref{420.8***}), we conclude that
\begin{eqnarray}\label{420.9***}
\begin{array}{ll}
\displaystyle \left| \Phi_{p_2} \left (\varphi (t_0), \cdots,
\varphi^{(p_2 -2)} (t_0), \varphi^{(p_2 -1)} (t_0)\right) -
\Phi_{p_2} \left (\varphi (1), \cdots, \varphi^{(p_2 -2)} (1),
\varphi^{(p_2 -1)} (1)\right) \right|
\\ \ns
\displaystyle = \left | \Phi_{p_2} \left (\varphi (t_0), \cdots,
\varphi^{(p_2 -2)} (t_0), a_{p_2}\right ) - \Phi_{p_2} \left (a_1,
\cdots, a_{p_2 -1}, 0\right ) \right |
\\ \ns
\displaystyle \ge \left | \Phi_{p_2} (a_1, \cdots, a_{p_2 -1},
a_{p_2}) - \Phi_{p_2} (a_1, \cdots, a_{p_2 -1}, 0) \right| \\ \ns
\displaystyle\quad - \left | \Phi_{p_2} \left (\varphi (t_0),
\cdots, \varphi^{(p_2 -2)} (t_0), a_{p_2}\right ) - \Phi_{p_2} (a_1,
\cdots, a_{p_2 -1}, a_{p_2}) \right| \\ \ns
\displaystyle \ge |b| - \frac{|b|}{2} = \frac{|b|}{2}.
\end{array}
\end{eqnarray}
By (\ref{418.5**}) and (\ref{421***}), one finds $\displaystyle
y_{q_1+p_2}(t)=\Phi_{p_2} \left(\varphi (t),
\overset{\cdot}{\varphi} (t), \cdots, \varphi^{(p_2 -1)}
(t)\right)$,  which is the $(q_1+p_2)$-th component of $y(t)$. Since
$u(t)=\left(0, \varphi^{(p_2)} (t)\right)$ is $C^{\infty}$ in $t $,
the solution $x(t)$ of (\ref{416***}) is also $C^{\infty}$. Hence
the corresponding control $v(t)=G(x(t),u(t))$ of (\ref{417***}) is
continuous, and therefore $\Phi_{p_2} \left(\varphi (t),
\overset{\cdot}{\varphi} (t), \cdots, \varphi^{(p_2 -1)} (t)\right)$
is differential with respect to $t $. By Differential Mean Value
Theorem, there is a $\xi \in (1,t_0)$ such that
\begin{eqnarray*}
\begin{array}{ll}
\displaystyle \left. \frac{d \Phi_{p_2} \left(\varphi (t),
\overset{\cdot}{\varphi}(t), \cdots, \varphi^{(p_2 -1)} (t)\right)}{
d t} \right|_{t=\xi} \\ \ns
\displaystyle = \frac{\Phi_{p_2} \left(\varphi(t_0),
\overset{\cdot}{\varphi} (t_0), \cdots, \varphi^{(p_2 -1)}
(t_0)\right) - \Phi_{p_2} \left (\varphi (1),
\overset{\cdot}{\varphi} (1), \cdots, \varphi^{(p_2 -1)} (1)
\right)} { t_0 -1 }.
\end{array}
\end{eqnarray*}
Combining the above formula and (\ref{420.9***}), recalling the
definition of $t_0$, we find
\begin{eqnarray} \label{424***}
\begin{array}{ll}
\displaystyle \left| \left. \frac{d \Phi_{p_2} \left(\varphi (t),
\overset{\cdot}{\varphi} (t), \cdots, \varphi^{(p_2 -1)}
(t)\right)}{ d t} \right|_{t=\xi} \right| \ge  \frac{|b|}{2(t_0 -1)}
\ge 2c > c.
\end{array}
\end{eqnarray}
By (\ref{425***}) and (\ref{424***}), it follows
$$\left |\left.\Phi_{p_2 +1} \left(\varphi (t),
 \cdots, \varphi^{(p_2 -2)} (t), \varphi^{(p_2 -1)} (t)\right)\right|_{t=\xi} \right |
>c,$$
which contradicts (\ref{420.7***}). Hence, (\ref{419***}) holds.

\medskip

{\it Step 3.} From the definition of  $\Phi_{p_2} (\lambda_1,
\cdots, \lambda_{p_2})$ and noting (\ref{419***}), we get
$$
h_{q_1 + p_2} (\overset{p_1 \text{ times}}{\overbrace{0, \cdots,
0}},\ \overset{p_2 \text{ times}}{\overbrace{ \lambda_1, \cdots,
\lambda_{p_2 -1}, \lambda_{p_2} }})\equiv h_{q_1 + p_2}
((\overset{p_1 \text{ times}}{\overbrace{0, \cdots, 0}},\
\overset{p_2 \text{ times}}{\overbrace{ \lambda_1, \cdots,
\lambda_{p_2 -1}, 0}}).
$$
Similarly,
\begin{eqnarray*}
\begin{array}{ll}
h_{q_1 + p_2 -1} (0, \cdots, 0, \lambda_1, \cdots, \lambda_{p_2 -1},
\lambda_{p_2}) \equiv h_{q_1 + p_2 -1} (0, \cdots, 0, \lambda_1,
\cdots, \lambda_{p_2 -1}, 0), \\ \ns \cdots \\ \ns
h_{q_1 +1} (0, \cdots, 0, \lambda_1, \cdots, \lambda_{p_2 -1},
\lambda_{p_2}) \equiv h_{q_1 + 1} (0, \cdots, 0, \lambda_1, \cdots,
\lambda_{p_2 -1}, 0).
\end{array}
\end{eqnarray*}
Thus, $y_{q_1 +1} (t),\  y_{q_1 +2} (t), \cdots, y_{q_1 + p_2} (t)$
only depend on $ x_{p_1 +1} (t),\  x_{p_1 +2} (t), \cdots, x_{p_1 +
p_2 -1} (t) $. Similar to the above, we deduce that
\begin{eqnarray*}
\begin{array}{ll}
h_{q_1 + p_2 -1} (0, \cdots, 0, \lambda_1, \cdots, \lambda_{p_2 -2},
\lambda_{p_2 -1}, \lambda_{p_2}) \equiv h_{q_1 + p_2} (0, \cdots, 0,
\lambda_1, \cdots, \lambda_{p_2 -2}, 0, 0), \\ \ns
 \cdots \\ \ns
h_{q_1 +1} (0, \cdots, 0, \lambda_1, \cdots, \lambda_{p_2 -2},
\lambda_{p_2 -1}, \lambda_{p_2}) \equiv h_{q_1 +1} (0, \cdots, 0,
\lambda_1, \cdots, \lambda_{p_2 -2}, 0, 0).
\end{array}
\end{eqnarray*}
Using the same method, we see that for $i= 0,1, \cdots, p_2 -2$, it
holds
\begin{eqnarray*}
\begin{array}{ll}
h_{q_1 + p_2 -i} (0, \cdots, 0, \lambda_1, \cdots, \lambda_{p_2 -i
-1}, \lambda_{p_2 -i}, \cdots, \lambda_{p_2 -1}, \lambda_{p_2}) \\
\ns
 \equiv  h_{q_1 + p_2 -i} (0, \cdots, 0, \lambda_1, \cdots,
\lambda_{p_2 -i -1}, 0, \cdots, 0),\\ \ns \cdots\\ \ns
h_{q_1 +1} (0, \cdots, 0, \lambda_1, \cdots, \lambda_{p_2 -i -1},
\lambda_{p_2 -i}, \cdots, \lambda_{p_2 -1}, \lambda_{p_2}) \\
\ns\equiv h_{q_1 +1} (0, \cdots, 0, \lambda_1, \cdots, \lambda_{p_2
-i -1}, 0, \cdots, 0).
\end{array}
\end{eqnarray*}

For $i= p_2 -2$, we have
$$
h_{q_1 +2} (0, \cdots, 0, \lambda_1, \cdots, \lambda_{p_2}) \equiv
h_{q_1 +2} (0, \cdots, 0, \lambda_1, 0, \cdots, 0),
$$
which leads to
\begin{eqnarray*}
\begin{array}{ll}
y_{q_1+2}(t) \equiv h_{q_1 +2} (0, \cdots, 0, x_{p_1 +1} (t), 0,
\cdots, 0).
\end{array}
\end{eqnarray*}
Then it follows from Lemma \ref{l41} that $y_{q_1 +2} (t) \equiv 0$.
By Lemma \ref{l42} and keeping in mind that $y(0)=H(x(0))=0$, we
obtain
\begin{eqnarray*}
\begin{array}{ll}
y_{q_1 +1} (t) =\displaystyle \int_0^t y_{q_1 +2} (\tau) d{\tau}
\equiv 0,\\ \ns
\ y_{q_1 +3} (t) = \overset{\cdot}{y}_{q_1+2} (t) \equiv 0, \
\cdots, \ y_{q_1 + q_2} (t) = y_{q_1 + 2}^{(q_2-2)} (t) \equiv 0,
\quad \forall t \ge 0.
\end{array}
\end{eqnarray*}

\medskip

{\it Step 4.} Similar to Step 2 and Step 3, we can prove that
$$
y_1 (t) = y_2 (t) = \cdots = y_{q_1} (t) \equiv  0, \quad \forall t
\ge 0.
$$
Therefore, $y(t) \equiv 0$. That is, the corresponding solution of
(\ref{417***}) is a fixed point in $\mathbb{R}^{q_1 + q_2}$
($=\mathbb{R}^{p_1 + p_2}$). But from (\ref{418***}), we see that
the solution of (\ref{416***}) is a curve in $\mathbb{R}^{p_1 +
p_2}$, which contradicts the fact that $H(x)$ is a homeomorphism
from $\mathbb{R}^{p_1 + p_2}$ to $\mathbb{R}^{p_1 + p_2}$.\endpf

\bigskip

{\it Proof of Proposition \ref{p45}.} It follows from Proposition
\ref{p43} that $r=\text{rank}D_1=\text{rank}D_2=s$. Denote $\{ p_i
\}_{i=1}^{m}=P(C_1,D_1)$ and $\{ q_i \}_{i=1}^{m}=P(C_2,D_2)$. By
Lemma \ref{l44}, it suffices to prove that $P(C_1,D_1)=P(C_2,D_2)$,
i.e.
\begin{equation} \label{430}
p_i = q_i, \quad  i=1, \cdots, r.
\end{equation}

We use the contradiction argument and suppose that (\ref{430}) is
false. Due to Lemma \ref{l21}, there exists a integer $\ell>0$ such
that
\begin{equation}\label{430.5}
0<p_r = q_r \le p_{r-1} = q_{r-1} \le \cdots \le
p_{r-\ell-1}=q_{r-\ell-1}, \quad p_{r-\ell} \neq q_{r-\ell}.
\end{equation}
Without loss of generality, we assume that $p_{r-\ell} <
q_{r-\ell}$.

Denote  by $(H(x), G(x,u))$ the topological equivalence
transformation from system
\begin{equation}\label{435}\overset{\cdot }{x}(t)=C_1x(t)+D_1u(t)
\end{equation} to system
\begin{equation}\label{436}\overset{\cdot }{y}(t)=C_2y(t)+D_2v(t).
\end{equation} Consider the following two cases: $p_{r-\ell}=1$ and $p_{r-\ell}\geq 2$.

\medskip

{\it Case 1.} $p_{r-\ell}=1$. Then (\ref{430.5}) leads to
$$
1= p_r = q_r = \cdots = p_{r-\ell-1}=q_{r-\ell-1} = p_{r-\ell} < 2
\le q_{r-\ell}.
$$

It is easy to see that the solutions of system (\ref{435})
associated to $x(0) = 0$ and\linebreak
$$u(t) = \Big(\overset{r-\ell-1 \text{ times}}{\overbrace{0, \cdots,
0}},\ \overset{\ell+1 \text{ times}}{\overbrace{ u_{r-\ell} (t),
u_{r-\ell+1} (t), \cdots, u_r (t)}},\  \overset{m-r \text{
times}}{\overbrace{0,\ \cdots,\ 0}}\Big)\in
C([0,+\infty);\mathbb{R}^m) $$ can be expressed as follows
\begin{eqnarray*}\label{431}
\begin{array}{ll}
x(t)&= \displaystyle \left(0, \cdots, 0, x_{k-\ell} (t),
x_{k-\ell+1} (t), \cdots, x_{k} (t)\right)^{\top} \\ \ns &=
\displaystyle \left(0, \cdots, 0, \int_0^t u_{r-\ell} (\tau) d\tau,
\int_0^t u_{r-\ell+1} (\tau) d\tau, \cdots, \int_0^t u_r(\tau) d\tau
\right)^{\top}.
\end{array}
\end{eqnarray*}
Denote by $H(x(t)) = y(t) = (y_1 (t), \cdots, y_k (t))^{\top}$  the
corresponding solution of (\ref{436}) by the transformation $H(x) =
(h_1(x), \cdots, h_k(x))^{\top}$.

\smallskip

Noting $q_1\geq q_2\geq\cdots\geq q_{r-\ell}\geq 2$, from Lemma
\ref{l41}, we deduce that
$$
y_{q_1} (t) = y_{q_1 + q_2} (t) = \cdots = y_{\sum_{j=1}^{r-\ell}
q_j} (t) \equiv 0, \quad \forall t  \ge 0.
$$
Therefore, by Lemma \ref{l42} and keeping in mind that
$y(0)=H(x(0))=0$, we obtain
\begin{equation} \label{432}
y_1 (t) = y_2 (t) = \cdots = y_{\sum_{j=1}^{r-\ell} q_j} (t) \equiv
0, \quad \forall t \ge 0.
\end{equation}

Since (\ref{435}) is completely controllable, for any fixed $t_0>0$,
we see that
$$\Theta_3=\{x(t_0;0,u(\cdot));\; u(t) = \left(0, \cdots,
0, u_{r-\ell} (t), u_{r-\ell+1} (t), \cdots, u_r (t), 0, \cdots,
0\right)\in C([0,+\infty);\mathbb{R}^m) \}$$ is a
$(\ell+1)$-dimensional topological manifold contained in $\mathbb
R^{k}$. On the other hand, (\ref{432}) implies that $y(t) =
H(\Theta_3)$ is at most a $(\sum_{j=r-\ell+1}^{r} q_j)$-dimensional
topological manifold contained in $\mathbb {R}^k$. Noting
$\sum_{j=r-\ell+1}^{r} q_j = \ell < \ell+1$, which contradicts the
fact that $H(x)$ is a homeomorphism.

\medskip

{\it Case 2.}  $p_{r-\ell} \ge 2$. For convenience, we assume
$p_{r-\ell} = 2$ and $p_r=1$ (for $p_{r-\ell} > 2$ or $p_r
> 1$, the proof is similar). Thus, (\ref{430.5}) gives
$$
1=p_r = q_r \le p_{r-1} = q_{r-1} \le \cdots \le
p_{r-\ell-1}=q_{r-\ell-1} \le p_{r-\ell}=2 <  q_{r-\ell}.
$$
The solutions of (\ref{435}) associated to $u(t) =
\Big(\overset{r-\ell-1 \text{ times}}{\overbrace{0, \cdots, 0}},\
\overset{\ell+1 \text{ times}}{\overbrace{ u_{r-\ell} (t),
u_{r-\ell+1} (t), \cdots, u_r (t)}},\  \overset{m-r \text{
times}}{\overbrace{0,\ \cdots,\ 0}}\Big)\linebreak\in
C([0,+\infty);\mathbb{R}^m)$ and $x(0) = 0$ is
\begin{eqnarray*}\label{433}
\begin{array}{ll}
x(t)&= \displaystyle \Big(0, \cdots, 0,\
\overset{\sum_{j=r-\ell}^{r} p_j \text{
times}}{\overbrace{x_{\sum_{j=1}^{r-\ell-1} p_j+1} (t),
x_{\sum_{j=1}^{r-\ell-1} p_j+2} (t), \cdots, x_{k} (t)}}\Big)^{\top}
\\ \ns &= \displaystyle \left(0, \cdots, 0, \int_0^t \int_0^{\tau}
u_{r-\ell} (\varsigma) d {\varsigma} d{\tau}, \int_0^t u_{r-\ell}
(\tau) d{\tau}, \cdots, \int_0^t u_r(\tau) d\tau \right)^{\top}.
\end{array}
\end{eqnarray*}

Similar to the proof of Lemma \ref{p44}, we arrive at
\begin{eqnarray*}
\begin{array}{ll}
 y_{\sum_{j=1}^{r-\ell-1} q_j+2} (t) \\ \ns =
h_{\sum_{j=1}^{r-\ell-1} q_j+2} (0, \cdots, 0,
x_{\sum_{j=1}^{r-\ell-1} p_j+1} (t), x_{\sum_{j=1}^{r-\ell-1} p_j+2}
(t), \cdots, x_{k} (t))
\\ \ns =  h_{\sum_{j=1}^{r-\ell-1} q_j+2} (0, \cdots, 0,
x_{\sum_{j=1}^{r-\ell-1} p_j+1} (t), 0, \cdots, x_{k} (t)),
\end{array}
\end{eqnarray*}
where $(0, \cdots, 0, x_{\sum_{j=1}^{r-\ell-1} p_j+1} (t), 0,
\cdots, x_{k} (t))^{\top}$ has the form of $\bar{x}$ in Lemma
\ref{l41}. Applying Lemma \ref{l41}, one finds
\begin{equation*}
y_{\sum_{j=1}^{r-\ell-1} q_j+2} (t) =h_{\sum_{j=1}^{r-\ell-1} q_j+2}
(0, \cdots, 0, x_{\sum_{j=1}^{r-\ell-1} p_j+1} (t), 0, \cdots, x_{k}
(t)) \equiv 0,\quad \forall t \ge 0.
\end{equation*}
Then by Lemma \ref{l42} and noting that $y(0)=H(x(0))=0$, we obtain
\begin{equation*}
y_{\sum_{j=1}^{r-\ell} q_j} (t) = y_{\sum_{j=1}^{r-\ell} q_j-1} (t)
= \cdots = y_{\sum_{j=1}^{r-\ell-1} q_j+1} (t) \equiv 0,\quad
\forall t \ge 0.
\end{equation*}
Similarly, for $i=1, \cdots, r-\ell-1$, it holds
\begin{equation*}
y_{\sum_{j=1}^{i} q_j} (t) = y_{\sum_{j=1}^{i} q_j-1} (t) =\cdots =
y_{\sum_{j=1}^{i-1} q_j+1} (t) \equiv 0,\quad \forall t \ge 0.
\end{equation*}
Hence, it is easy to see that
\begin{equation} \label{434}
y_1 (t) = y_2 (t) = \cdots = y_{\sum_{j=1}^{r-\ell} q_j-1} (t) =
y_{\sum_{j=1}^{r-\ell} q_j} (t) \equiv 0, \quad \forall t \ge 0.
\end{equation}

Since (\ref{435}) is completely controllable,  for any fixed
$t_0>0$, we see that
$$\Theta_4=\{x(t_0;0,u(\cdot));\; u(t) = \left(0, \cdots,
0, u_{r-\ell} (t), u_{r-\ell+1} (t), \cdots, u_r (t), 0, \cdots,
0\right)\in C([0,+\infty);\mathbb{R}^m) \}$$ is a
$(\sum_{j=r-\ell}^{r} p_j)$-dimensional topological manifold
contained in $\mathbb {R}^{k}$. On the other hand, (\ref{434})
implies that $y(t) = H(\Theta_4)$ is at most a
$(\sum_{j=r-\ell+1}^{r} q_j)$-dimensional topological manifold
contained in $\mathbb {R}^k$. Noting $\sum_{j=r-\ell+1}^{r} q_j =
\sum_{j=r-\ell+1}^{r} p_j < \sum_{j=r-\ell}^{r} p_j$, which
contradicts the fact that $H(x)$ is a homeomorphism.
\endpf

\section{Appendix A: Proof of Proposition \ref{p11}}

This section is devoted to a proof of Proposition \ref{p11}. The
proof is divided into several steps.

\medskip

{\it Step 1.} We use the contradiction argument to prove that
$H(x,u)$ is independent of $u$. Assume that there exist a $x^0\in
\mathbb{R}^n$ and a $u^0\in\mathbb{R}^m\backslash\{0\}$ such that
$$H\left(x^0, u^0\right) \neq H\left(x^0, 0\right).$$

Denote $\varepsilon_0 = \left| H\left(x^0, u^0\right) - H\left(x^0,
0\right) \right|$. It suffices to find a $x^1\in \mathbb{R}^n$
satisfying the following two estimates (the existence of $x^1$ will
be proved in Step 2)
\begin{eqnarray}\label{14}\left| H\left(x^1, u^0\right) - H\left(x^0, u^0\right) \right| <
\frac{\varepsilon_0}{3},\end{eqnarray}
\begin{eqnarray}\label{15}\left| H\left(x^1, u^0\right) - H\left(x^0, 0\right) \right| <
\frac{\varepsilon_0}{3}.\end{eqnarray}
Indeed, by (\ref{14}) and (\ref{15}), it follows
\begin{eqnarray*}\begin{array}{ll}\displaystyle\left| H\left(x^0,
u^0\right) - H\left(x^0, 0\right) \right|\\ \ns \displaystyle\leq
\left| H\left(x^0, u^0\right) - H\left(x^1, u^0\right) \right| +
\left| H\left(x^1, u^0\right) - H\left(x^0, 0\right) \right|\\ \ns <
\displaystyle \frac{2 \varepsilon_0}{3},\end{array}\end{eqnarray*}
which contradicts the fact that $\left| H\left(x^0, u^0\right) -
H\left(x^0, 0\right) \right|=\varepsilon_0$.

\medskip

{\it Step 2.} We now show that there exists a $x^1\in \mathbb{R}^n$
such that (\ref{14}) and ({\ref{15}) hold simultaneously.

Recalling that $H(x,u)$ is continuous in $x$, we conclude that there
exists a constant $\delta>0$ such that for any $\left| x - x^0
\right| < \delta$, it holds
\begin{equation}\label{110}\left| H\left(x,
u^0\right) - H\left(x^0, u^0\right) \right| <
\frac{\varepsilon_0}{3}.\end{equation}
Construct a control as follows
\begin{eqnarray*}u(t)=\left\{\begin{array}{ll}
\displaystyle \frac{tu^0}{t_0},\quad 0\leq t\leq t_0,\\ \ns u^0,\ \
\quad t>t_0,
\end{array}\right.
\end{eqnarray*}
where $t_0>0$ will be given later.
Clearly,
\begin{eqnarray}\label{16}\left|u(t)\right|\leq\left|u^0\right|,\quad
\forall t\geq 0.\end{eqnarray}

The solution of system (\ref{12}) associated to $x(0)=x^0$ and
$u(t)$ is expressed as follows
\begin{eqnarray}\label{115}
 x(t)=e^{A_1t}x^0+\int_0^t e^{A_1(t-s)}B_1u(s)ds.
\end{eqnarray}
By (\ref{16}), one can find a $t_1 > 0$ (which is independent of
$t_0$) small enough such that
\begin{eqnarray}\label{17}\left| x(t) - x^0 \right| < \delta,\quad
\forall 0\leq t \leq t_1.\end{eqnarray}
By (\ref{110}) and (\ref{17}), we obtain
\begin{eqnarray}\label{113}\left| H\left(x(t), u^0\right) - H\left(x^0, u^0\right) \right| <
\frac{\varepsilon_0}{3},\quad \forall 0\leq t \leq
t_1.\end{eqnarray}

On the other hand, from (\ref{17}), we get
\begin{eqnarray}\label{18}\left| x(t) \right| <  \left| x^0 \right| + \delta,\quad
\forall 0\leq t \leq t_1.\end{eqnarray}
Since $G(x,u)$ maps a bounded set in
$\mathbb{R}^n\times\mathbb{R}^m$ to a bounded set in
$\mathbb{R}^n\times\mathbb{R}^m$, (\ref{16}) and (\ref{18}) imply
that there exists a constant $c>0$ satisfying
\begin{equation}\label{111}\left| v(t) \right|= \left| G(x(t),u(t)) \right|< c,\quad
\forall 0\leq t \leq t_1.\end{equation}
The solution of system (\ref{13}) corresponding to $x(t)$ and $u(t)$
is expressed as follows
\begin{eqnarray}\label{19}
 H(x(t),u(t))=y(t)=e^{A_2t}y^0+\int_0^t e^{A_2(t-s)}B_2v(s)ds.
\end{eqnarray}
Taking $t=0$ in (\ref{19}), we see that $
H\left(x^0,u(0)\right)=y^0. $ Recalling $u(0)=0$, (\ref{19}) can be
rewritten as:
$$
 H(x(t), u(t)) = e^{A_2t}H\left(x^0, 0\right) +\int_0^t e^{A_2(t-s)}B_2v(s)ds.
$$
The above formula and (\ref{111}) imply that there exists a
$t_2\in(0, t_1)$ (which is independent of $t_0$) such that
\begin{equation}\label{114}\left|
H(x(t), u(t))- H\left(x^0, 0\right) \right| <
\frac{\varepsilon_0}{3},\quad \forall 0\leq t\leq t_2.\end{equation}

Setting $t_0=t_2/2$, from the construction of $u(t)$ we have
$u(t_2)=u^0$. Taking $t = t_2$, $x^1 = x(t_2)$ in (\ref{113}) and
(\ref{114}), we obtain that inequalities (\ref{14}) and (\ref{15})
hold simultaneously.

\medskip

{\it Step 3.} Since $H(x,u)$ is independent of $u$, we simply denote
it as $H(x)$. We show that $H(x)$ is a homeomorphism from
$\mathbb{R}^n$ to $\mathbb{R}^n$.

Recall the inverse transformation of $F(x,u)$ is
$F^{-1}(y,v)=(Z(y,v),W(y,v))$. Similar to Step 1, we see that
$$Z(y,v)\equiv Z(y),\quad  \forall y\in\mathbb{R}^n,\ \forall
v\in\mathbb{R}^m.$$

We claim that $H(x)$ is invertible and $H^{-1}=Z$. In fact, since
$F(\cdot,\cdot)$ is a homeomorphism from
$\mathbb{R}^n\times\mathbb{R}^m$ to
$\mathbb{R}^n\times\mathbb{R}^m$, we obtain
$$(x,u)=F^{-1}(F(x,u))=(Z(H(x)),W(H(x),G(x,u))),\quad \forall x\in\mathbb{R}^n,\ \forall u\in\mathbb{R}^m,$$
which yields $$x=Z(H(x)),\quad \forall x\in\mathbb{R}^n.$$
Similarly, we have
$$y=H(Z(y)),\quad \forall y\in\mathbb{R}^n.$$ Hence, $H^{-1}=Z$.

Furthermore, noting $H$ and $Z$ are continuous functions, we
conclude that $H(x)$ is a homeomorphism from $\mathbb{R}^n$ to
$\mathbb{R}^n$.
\endpf

\begin{remark}\label{r15}
Proposition \ref{p11} shows that the state variable and the control
variable play different roles in the topological transformation
function. This is mainly due to the formulae (\ref{115}) and
(\ref{19}), which show the relations among the variables $x$, $u$
and $y$. More precisely, if $u$ change a lot in a very short time
duration, by (\ref{115}) (resp. (\ref{19})), $x$ (resp. $y$) might
change little. This inspire us to guess that $y$ do not depend on
$u$.
\end{remark}

\begin{remark}\label{r16}
To prove Proposition \ref{p11}, we only need that $G(x,u)$ maps a
bounded set in $\mathbb{R}^n\times\mathbb{R}^m$ to a bounded set in
$\mathbb{R}^n\times\mathbb{R}^m$ and $H(x,u)$, $Z(y,v)$ are
continuous on the first variable.
\end{remark}

\end {document}